\documentclass{amsart}

\usepackage{mathrsfs}
\usepackage[all]{xy}
\usepackage{amssymb,latexsym,amsmath,amsthm, mathrsfs}
\usepackage{graphicx}
\usepackage{setspace}
\newtheorem{theorem}{Theorem}[section]
\newtheorem{lemma}[theorem]{Lemma}
\newtheorem{definition}[theorem]{Definition}
\newtheorem{example}[theorem]{Example}

\newtheorem{proposition}[theorem]{Proposition}
\newtheorem{corollary}[theorem]{Corollary}
\newenvironment{acknowledgement}[1][Acknowledgements]
{\begin{trivlist} \item[\hskip \labelsep {\bfseries #1}]}
{\end{trivlist}}

\DeclareMathSymbol{\N}{\mathbin}{AMSb}{"4E}
\DeclareMathSymbol{\Z}{\mathbin}{AMSb}{"5A}
\DeclareMathSymbol{\R}{\mathbin}{AMSb}{"52}
\DeclareMathSymbol{\Q}{\mathbin}{AMSb}{"51}
\DeclareMathSymbol{\I}{\mathbin}{AMSb}{"49}
\DeclareMathSymbol{\C}{\mathbin}{AMSb}{"43}

\def\P{{\mathbb P}}
\def\RP{{\mathbb R \mathbb P}}

\def\A{{\mathbb A}}

\numberwithin{equation}{section}

\begin{document}

\title{Harnack-Thom Theorem for higher cycle groups and Picard varieties}
\author{Jyh-Haur Teh}

\address{Department of Mathematics, National Tsing Hua University of Taiwan, No. 101, Kuang Fu Road,
Hsinchu, 30043, Taiwan.}

\email{jyhhaur@math.nthu.edu.tw}


\date{}


\begin{abstract}
We generalize the Harnack-Thom theorem to relate the ranks of the
Lawson homology groups with $\Z_2$-coefficients of a real
quasiprojective variety with the ranks of its reduced real Lawson
homology groups. In the case of zero-cycle group, we recover the
classical Harnack-Thom theorem and generalize the classical version
to include real quasiprojective varieties. We use Weil's
construction of Picard varieties to construct reduced real Picard
groups, and Milnor's construction of universal bundles to construct
some weak models of classifying spaces of some cycle groups. These
weak models are used to produce long exact sequences of homotopy
groups which are the main tool in computing the homotopy groups of
some cycle groups of divisors. We obtain some congruences involving
the Picard number of a nonsingular real projective variety and the
rank of its reduced real Lawson homology groups of divisors.
\end{abstract}

\maketitle

\section{Introduction}
In \cite{F1, FL1, L1}, Friedlander and Lawson constructed Lawson
homology and morphic cohomology, which serve as an enrichment of
singular homology and singular cohomology, respectively, for complex
projective varieties. In \cite{T1}, the author constructed parallel
theories for real projective varieties which are called reduced real
Lawson homology and reduced real morphic cohomology. They enjoy many
nice properties such as the Lawson suspension property, the homotopy
invariance property, the bundle projection property, the splitting
properties and for each theory there exists a localization long
exact sequence. By using the Friedlander-Lawson moving lemma (see
\cite{FL2}), it is shown that there is a duality theorem between
Lawson homology and morphic cohomology (see \cite{FL3}) and a
duality theorem between reduced real Lawson homology and reduced
real morphic cohomology (see \cite{T1}). Furthermore, this duality
is compatible with Poincar\'{e} duality.

The Harnack theorem says that a nonsingular totally real curve of
degree $d$ in $\RP^2$ has at most $g(d)+1$ connected components
where $g(d)=\frac{(d-1)(d-2)}{2}$. Later on Thom generalized
Harnack's result to a statement which says that for a real
projective variety $X$, the total Betti number $B(X), B(ReX)$ and
the Euler characteristic $\chi(X), \chi(ReX)$ in $\Z_2$-coefficients
of $X$ and the real points $ReX$ of $X$ respectively satisfy the
following relations(see \cite{DK, Gud, Thom}):
$$B(ReX)\leq B(X), B(ReX)\equiv B(X) \mbox{mod }2,
\chi(ReX)\equiv \chi(X) \mbox{mod }2$$

In section 2 we give an overview of Lawson homology and reduced real
Lawson homology. In section 3 we prove a splitting theorem which is
the core of the proof of our main theorem. In section 4 we extend
the classical Harnack-Thom theorem to a statement involving the
ranks of Lawson homology groups with $\Z_2$-coefficients and the
ranks of reduced real Lawson homology groups. For 0-cycle groups, we
recover the Harnack-Thom theorem and generalize it to real
quasiprojective varieties, in which case we need to use Borel-Moore
homology instead of singular homology. To construct some nontrivial
examples, we apply Weil's construction of Picard varieties to
construct reduced real Picard groups in section 5. In section 6 we
prove a vanishing theorem for the reduced real Lawson homology
groups of divisors and under some mild conditions, we get the
following result by applying our main theorem from section 4:
$$\rho(X)+1\equiv
dim_{\Z_2}\pi_0R_{m-1}(X)+dim_{\Z_2}\pi_1R_{m-1}(X) mod 2$$ where
$\rho(X)$ is the Picard number of $X$.

The results of this paper suggest that Lawson homology and reduced
real Lawson homology are useful enrichments of singular homology.

\section{Review of Lawson homology and reduced real Lawson
homology} Let us recall some basic properties of Lawson homology and
reduced real Lawson homology (see \cite{F1, L1, T1}). For a
projective variety $X$, denote the set of effective $p$-cycles of
degree $d$ by $\mathscr{C}_{p, d}(X)$. By the Chow Theorem (see
\cite{S}), $\mathscr{C}_{p, d}(X)$ can be realized as a complex
projective variety. With the analytic topology on $\mathscr{C}_{p,
d}(X)$, we get a compact topological space $K_{p,
d}(X)=\coprod_{d_1+d_2\leq d}\mathscr{C}_{p, d_1}(X)\times
\mathscr{C}_{p, d_2}(X)/\sim$ where $\sim$ is the equivalence
relation defined by $(a, b)\sim (c, d)$ if and only if $a+d=b+c$.
These spaces form a filtration:
$$K_{p, 0}(X)\subset K_{p, 1}(X) \subset K_{p, 2}(X) \subset
\cdots  =Z_p(X)$$ where $Z_p(X)$ is the naive group completion of
the monoid $\mathscr{C}_p(X)=\coprod_{d\geq 0}\mathscr{C}_{p,
d}(X)$. We give $Z_p(X)$ the weak topology defined by this
filtration, i.e., $U\subset Z_p(X)$ is open if and only if $U\cap
K_{p, d}(X)$ is open for all $d$. We define the $n$-th Lawson
homology group of $p$-cycles to be
$$L_pH_n(X)=\pi_{n-2p}Z_p(X)$$ the $(n-2p)$-th homotopy group of $Z_p(X)$.
We define the $n$-th Lawson homology group with $\Z_2$-coefficients
of $p$-cycles to be
$$L_pH_n(X; \Z_2)=\pi_{n-2p}\left(\frac{Z_p(X)}{2Z_p(X)}\right).$$
For the zero-cycle group, by the Dold-Thom Theorem, we have an
isomorphism
$$L_0H_n(X; \Z_2)=H_n(X; \Z_2)$$ between Lawson homology and
singular homology.

For a quasiprojective variety $U$, there exist projective varieties
$X$ and $Y$ where $Y\subset X$ such that $U=X-Y$. The Lawson
homology group of $U$ is defined to be
$$L_pH_n(U)=\pi_{n-2p}\left(\frac{Z_p(X)}{Z_p(Y)}\right).$$
It is proved in \cite{Li1} that this definition is independent of
the choice of $X$ and $Y$. We define Lawson homology with
$\Z_2$-coefficients of $U$ to be
$$L_pH_n(U;
\Z_2)=\pi_{n-2p}\left(\frac{\frac{Z_p(X)}{2Z_p(X)}}{\frac{Z_p(Y)}{2Z_p(Y)}}\right)=
\pi_{n-2p}\left(\frac{\frac{Z_p(X)}{Z_p(Y)}}{2\frac{Z_p(X)}{Z_p(Y)}}\right)=
\pi_{n-2p}\left(\frac{\frac{Z_p(X)}{Z_p(Y)}}{\frac{2Z_p(X)}{2Z_p(Y)}}\right)$$
For the zero-cycle group, by the Dold-Thom Theorem (see \cite{F3}
Proposition 1.6), we have an isomorphism
$$L_0H_n(U)=H^{BM}_n(U) \mbox{ and } L_0H_n(U; \Z_2)=H^{BM}_n(U; \Z_2)$$
where $H^{BM}_*$ denotes Borel-Moore homology.

A real projective variety $X\subset \P^n$ is a complex projective
variety which is invariant under conjugation. Equivalently, it is a
complex projective variety defined by some real polynomials.
Conjugation induces a $\Z_2$-action on $Z_p(X)$. Let $Z_p(X)_{\R}$
be the subgroup of $p$-cycles on $X$ which are invariant under this
action and let $Z_p(X)^{av}$ be the subgroup consisting of cycles of
the form $c+\overline{c}$ where $c\in Z_p(X)$ and $\overline{c}$ is
the conjugate cycle of $c$. These two subgroups are endowed with the
subspace topology. Define the reduced real $p$-cycle group to be
$$R_p(X)=\frac{Z_p(X)_{\R}}{Z_p(X)^{av}}.$$ It is shown in
\cite{T1} Proposition 2.3 that $Z_p(X)^{av}$ is a closed subgroup of
$Z_p(X)_{\R}$ and in the Appendix of \cite{T1} that all these cycle
groups $Z_p(X), Z_p(X)_{\R}, Z_p(X)^{av}, R_p(X)$ are CW-complexes.
We define the $n$-th reduced real Lawson homology group of
$p$-cycles to be
$$RL_pH_n(X)=\pi_{n-p}R_p(X).$$
For zero cycles, it is shown in \cite{T1} Proposition 2.7 that
$RL_0H_n(X)=H_n(ReX; \Z_2)$, the singular homology group of the real
points of $X$.

We define the reduced real Lawson homology group of a real
quasiprojective variety $U=X-Y$ to be
$$RL_pH_n(U)=\pi_{n-p}\left(\frac{R_p(X)}{R_p(Y)}\right)$$ where
$X, Y$ are real projective varieties and $Y\subset X$. It is shown
in \cite{T1} that this definition is independent of the choice of
$X, Y$. For the group of zero cycles, we have
$RL_oH_n(U)=\pi_n\left(\frac{R_0(X)}{R_0(Y)}\right)=H^{BM}_n(Re U;
\Z_2).$

\section{The Splitting Theorem}
Let us recall that the real part $RP(C)$ of a cycle $C$, roughly
speaking, is the part consisting of irreducible real subvarieties
and the averaged part $AP(C)$ of a cycle $C$ is the part consisting
of conjugate pairs of complex cycles. The imaginary part is the part
left after canceling out the real and averaged parts. We give the
precise definition in the following:
\begin{definition}
For any $f\in Z_p(X)$, let $f=\sum_{i\in I}n_iV_i$ be in the reduced
form, i.e., each $V_i$ is an irreducible subvariety of $X$ and
$V_i=V_j$ if and only if $i=j$. Let
$$RP(f)=\sum_{i \in I, \overline{V_i}=V_i}n_i V_i$$ which is
called the  real part of $f$. Let $$J=\{i \in I|V_i \mbox{ is not
real and } \overline{V_i} \mbox{ is also a component of } f \}$$ and
for $i\in J$, let $m_i$ be the maximum value of the coefficients of
$V_i$ and $\overline{V_i}$. Define the averaged part to be
$$AP(f)=\sum_{i\in J}m_i(V_i+\overline{V}_i)$$ and the
imaginary part to be
$$IP(f)=f-RP(f)-AP(f).$$
\label{real part}
\end{definition}
It is easy to see that $f$ is a real cycle if and only if $IP(f)=0$
and a real cycle $g$ is an averaged cycle if and only if $RP(g)$ is
divisible by 2.

In the following, we will assume that $X$ is a real projective
variety.

\begin{proposition}
The following sequence is exact: $0\longrightarrow
\frac{Z_p(X)_{\R}}{2Z_p(X)_{\R}} \overset{i}{\longrightarrow}
\frac{Z_p(X)}{2Z_p(X)} \overset{1+c_*}{\longrightarrow}
\frac{Z_p(X)^{av}}{2Z_p(X)_{\R}} \longrightarrow 0.$ \label{smith
exact sequence}
\end{proposition}

\begin{proof}
It is easy to check that $i(f+2Z_p(X)_{\R})=f+2Z_p(X)$ is well
defined and injective and
$(1+c_*)(f+2Z_p(X))=f+\bar{f}+2Z_p(X)_{\R}$ is well defined and
surjective. The map $1+c_*$ sends the image of $i$ to 0, thus the
only thing we need to prove is for $f+\bar{f}\in 2Z_p(X)_{\R}$, $f$
is in $Z_p(X)_{\R}$. Since
$f+\bar{f}=2RP(f)+2AP(f)+IP(f)+IP(\bar{f})\in 2Z_p(X)_{\R}$, this
implies that $IP(f)=IP(\bar{f})=0$ so $f\in Z_p(X)_{\R}$.
\end{proof}

\begin{definition}
Let \index{$Q_p(X)$} $Q_p(X)$ be the collection of all averaged
cycles $c$ such that there exists a sequence $\{v_i\}\subset
Z_p(X)_{\R}$ where $v_i=RP(v_i)$ for all $i$ and ${v_i}$ converges
to $c$. It is not difficult to see that $Q_p(X)$ is a topological
subgroup of $Z_p(X)^{av}$. Let \index{$ZQ_p(X)_{\R}$}
$ZQ_p(X)_{\R}=2Z_p(X)_{\R}+Q_p(X)$ denote the internal sum of
$2Z_p(X)_{\R}$ and $Q_p(X)$. Then $ZQ_p(X)_{\R}$ is again a
topological subgroup of $Z_p(X)^{av}$. The group $Q_p(X)$ is the
intersection of the closure of the group formed by irreducible real
$p$-subvarieties with the averaged $p$-cycle group. Thus
$ZQ_p(X)_{\R}$ is a closed subgroup.
\end{definition}

\begin{proposition}
For a real projective variety $X$, $ZQ_0(X)_{\R}=2Z_0(X)_{\R}$.
\label{o cycle}
\end{proposition}

\begin{proof}
The free abelian group $Z_0(ReX)$ generated by real points of $X$ is
closed in $Z_0(X)$, so if $c\in Q_0(X)$, then $c\in 2Z_0(X)_{\R}$
(see Proposition 2.7 in \cite{T1}).
\end{proof}

The following example was given by Lawson to show that the set of
1-cycles formed by irreducible real subvarieties may not be closed
which contrasts to the case of 0-cycles, i.e., $ZQ_p(X)_{\R}$ may
not equal to $2Z_p(X)_{\R}$ if $p>0$.

\begin{example}
In $\P^2$, consider the sequence of irreducible real subvarieties
$V_{\epsilon}=\mbox{ zero locus of } X^2+Y^2-\epsilon Z^2$. As
$\epsilon$ converges to 0, $V_{\epsilon}$ converges to the cycle
formed by two lines $X=iY$ and $X=-iY$ which is an averaged cycle
but not in $2Z_p(X)_{\R}$.
\end{example}

\begin{lemma}
Suppose that the sequence $\{AP(f_i)\}$ converges to $f$ where
$f_i\in Z_p(X)$, then $RP(f)\in 2Z_p(X)_{\R}$. \label{real part
converges}
\end{lemma}

\begin{proof}
Since $AP(f_i)\in Z_p(X)^{av}$ and $Z_p(X)^{av}$ is closed, $f$ is
in $Z_p(X)^{av}$ so $RP(f)\in 2Z_p(X)_{\R}$.
\end{proof}

\begin{lemma}
Define $\widetilde{AP}: \frac{Z_p(X)_{\R}}{ZQ_p(X)_{\R}}
\longrightarrow \frac{Z_p(X)^{av}}{ZQ_p(X)_{\R}}$ by
$$f+ZQ_p(X)_{\R}\longmapsto AP(f)+ZQ_p(X)_{\R}.$$ Then $\widetilde{AP}$ is continuous.
\end{lemma}

\begin{proof}
There is a filtration called the canonical real filtration
$$K_1 \subset K_2 \subset K_3 \subset \cdots =Z_p(X)_{\R}$$ where each $K_i$
is compact and the topology of $Z_p(X)_{\R}$ is given by the weak
topology induced from this filtration. Thus the filtration

$$K_1+ZQ_p(X)_{\R} \subset K_2+ZQ_p(X)_{\R} \subset \cdots
=\frac{Z_p(X)_{\R}}{ZQ_p(X)_{\R}}$$ defines the topology of
$\frac{Z_p(X)_{\R}}{ZQ_p(X)_{\R}}$.

To show that $\widetilde{AP}$ is continuous, it suffices to show
that $\widetilde{AP}$ maps convergent sequences to convergent
sequences. The reason is as follows: by elementary point-set
topology, $\widetilde{AP}$ is continuous if and only if
$\widetilde{AP}(\overline{A})\subset \overline{\widetilde{AP}(A)}$
for any subset $A\subset \frac{Z_p(X)_{\R}}{ZQ_p(X)_{\R}}$ where
$\overline{A}$ is the closure of $A$. Let
$\pi:Z_p(X)_{\R}\rightarrow \frac{Z_p(X)_{\R}}{ZQ_p(X)_{\R}}$ be the
quotient map. For any subset $A\subset
\frac{Z_p(X)_{\R}}{ZQ_p(X)_{\R}}$, since $\pi$ is an open map, we
have $\pi^{-1}(\overline{A})\subset \overline{\pi^{-1}(A)}$ where
$\overline{A}$ is the closure of $A$. Hence for $x\in \overline{A}$,
there is $y\in \overline{\pi^{-1}(A)}$ such that $\pi(y)=x$. Since
$Z_p(X)_{\R}$ is a CW-complex (actually it is also a metric space),
there is $y_n\in \pi^{-1}(A)$ such that $y_n\rightarrow y$. If
$\widetilde{AP}$ maps convergent sequences to convergent sequences,
we see that $\widetilde{AP}(\pi(y_n))$ converges to
$\widetilde{AP}(x)$. Since $\pi(y_n)\in A$ and $\pi(y)=x$, we have
$\widetilde{AP}(x)\in \overline{\widetilde{AP}(A)}$ which implies
that $\widetilde{AP}$ is continuous.

Suppose that $f_i+ZQ_p(X)_{\R}$ converges to $ZQ_p(X)_{\R}$. Since
$A=\{f_i+ZQ_p(X)_{\R}\}\cup \{ZQ_p(X)_{\R}\}$ is compact and
$\frac{Z_p(X)_{\R}}{ZQ_p(X)_{\R}}$ is Hausdorff, by Lemma 2.2 in
\cite{T1}, $A\subset K_n+ZQ_p(X)_{\R}$ for some $n$. Thus there
exists $g_i\in K_n$ such that under the quotient map $q$,
$q(g_i)=f_i+ZQ_p(X)_{\R}$ for all $i$. The set $K_n$ is compact,
thus $\{g_i\}$ has a convergent subsequence.

Let $\{g_{ij}\}$ be a subsequence of $\{g_i\}$ which converges to
$g$. Since $g_{ij}+ZQ_p(X)_{\R}$ converges to $ZQ_p(X)_{\R}$, we
have $g\in ZQ_p(X)_{\R}$. The set $\{g_{ij}\}\subset K_n$ and each
$g_{ij}$ is a real cycle which implies that $\{AP(g_{ij})\}\subset
K_n$ and hence $\{AP(g_{ij})\}$ has a convergent subsequence. Let
$\{AP(g_{ijk})\}$ be a subsequence of $\{AP(g_{ij})\}$ which
converges to a real cycle $h$. Since $\{g_{ijk}\}$ is a subsequence
of $\{g_{ij}\}$, it converges to $g$, hence
$$RP(g_{ijk})=g_{ijk}-AP(g_{ijk})\longrightarrow g-h.$$
By Lemma \ref{real part converges}, $RP(h)\in 2Z_p(X)_{\R}$ and
since $RP(g)\in 2Z_p(X)_{\R}$ we have $RP(g-h)\in 2Z_p(X)_{\R}$. The
cycle $g-h$ is a real cycle and this implies $g-h\in Z_p(X)^{av}$.
Furthermore, since $\{RP(g_{ijk})\}\longrightarrow g-h$, by
definition, $g-h\in Q_p(X)\subset ZQ_p(X)_{\R}$. The cycle $g$ is in
$ZQ_p(X)_{\R}$ thus $h\in ZQ_p(X)_{\R}$. Passing to the quotient, we
see that $AP(g_{ijk})+ZQ_p(X)_{\R}\longrightarrow
h+ZQ_p(X)_{\R}=ZQ_p(X)_{\R}$ for any convergent subsequence
$\{AP(g_{ijk})\}$ of $\{AP(g_{ij})\}$, thus $AP(g_{ij})+ZQ_p(X)_{\R}
\longrightarrow ZQ_p(X)_{\R}$. For any convergent subsequence
$\{g_{ij}\}$ of $\{g_i\}$,
$\widetilde{AP}(g_{ij}+ZQ_p(X)_{\R})=AP(g_{ij})+ZQ_p(X)_{\R}$
converges to the point $ZQ_p(X)_{\R}$. Consequently, this implies
that
$\widetilde{AP}(g_i+ZQ_p(X)_{\R})=\widetilde{AP}(f_i+ZQ_p(X)_{\R})$
converges to $ZQ_p(X)_{\R}$. So $\widetilde{AP}$ is continuous.
\end{proof}

\begin{lemma}
Define $\widetilde{RP}: \frac{Z_p(X)_{\R}}{Z_p(X)^{av}}
\longrightarrow \frac{Z_p(X)_{\R}}{ZQ_p(X)_{\R}}$ by
$$f+Z_p(X)^{av}\longmapsto RP(f)+ZQ_p(X)_{\R}.$$ Then $\widetilde{RP}$
is continuous.
\end{lemma}

\begin{proof}
We proceed as in the proof above. The canonical real filtration
$$K_1 \subset K_2 \subset K_3 \subset \cdots =Z_p(X)_{\R}$$ induces a filtration
$$K_1+Z_p(X)^{av} \subset K_2+Z_p(X)^{av} \subset \cdots
=\frac{Z_p(X)_{\R}}{Z_p(X)^{av}}$$ which defines the topology of
$\frac{Z_p(X)_{\R}}{Z_p(X)^{av}}$ and the filtration

$$K_1+ZQ_p(X)_{\R} \subset K_2+ZQ_p(X)_{\R} \subset \cdots
=\frac{Z_p(X)_{\R}}{ZQ_p(X)_{\R}}$$ which defines the topology of
$\frac{Z_p(X)_{\R}}{ZQ_p(X)_{\R}}$.

By an argument similar to that of the previous proof, it suffices to
prove that $\widetilde{RP}$ maps convergent sequences to convergent
sequences. Suppose that $f_i+Z_p(X)^{av}$ converges to
$Z_p(X)^{av}$. Since $A=\{f_i+Z_p(X)^{av}\}\cup \{Z_p(X)^{av}\}$ is
compact, we have $A\subset K_n+Z_p(X)^{av}$ for some $n$. Thus there
exists $g_i\in K_n$ such that under the quotient map $q$,
$q(g_i)=f_i+Z_p(X)^{av}$ for all $i$. Let $\{g_{ij}\}$ be a
subsequence of $\{g_i\}$ which converges to $g$. Since
$g_{ij}+Z_p(X)^{av}$ converges to $Z_p(X)^{av}$, this implies that
$g\in Z_p(X)^{av}$. The set $\{g_{ij}\}\subset K_n$ implies that
$\{RP(g_{ij})\}\subset K_n$ thus $\{RP(g_{ij})\}$ has a convergent
subsequence. Let $\{RP(g_{ijk})\}$ be a subsequence of
$\{RP(g_{ij})\}$ which converges to a real cycle $h$. Since
$\{g_{ijk}\}$ is a subsequence of $\{g_{ij}\}$, it converges to $g$,
hence
$$AP(g_{ijk})=g_{ijk}-RP(g_{ijk})\longrightarrow g-h.$$
By Lemma \ref{real part converges}, $RP(g-h)\in 2Z_p(X)_{\R}$ and
since $RP(g)\in 2Z_p(X)_{\R}$ we have $RP(h)\in 2Z_p(X)_{\R}$. The
cycle $h$ is a real cycle and thus $h\in Z_p(X)^{av}$. Furthermore,
$\{RP(g_{ijk})\}\longrightarrow h$, so by definition, $h\in Q_p(X)$.
Passing to the quotient, we see that
$RP(g_{ijk})+ZQ_p(X)_{\R}\longrightarrow
h+ZQ_p(X)_{\R}=ZQ_p(X)_{\R}$. Thus $RP(g_{ij})+ZQ_p(X)_{\R}
\longrightarrow ZQ_p(X)_{\R}$. For every convergent subsequence
$\{g_{ij}\}$ of $\{g_i\}$, $\widetilde{RP}(g_{ij}+Z_p(X)^{av})$
converges to the point $ZQ_p(X)_{\R}$, thus
$\widetilde{RP}(g_i+Z_p(X)^{av})\longrightarrow ZQ_p(X)_{\R}$ and
therefore $\widetilde{RP}$ is continuous.
\end{proof}

\begin{theorem}(The splitting theorem)
$\frac{Z_p(X)_{\R}}{ZQ_p(X)_{\R}}$ is isomorphic as a topological
group to $\frac{Z_p(X)_{\R}}{Z_p(X)^{av}}\times
\frac{Z_p(X)^{av}}{ZQ_p(X)_{\R}}$ \label{split}
\end{theorem}

\begin{proof}
Define $\psi: \frac{Z_p(X)_{\R}}{ZQ_p(X)_{\R}} \longrightarrow
\frac{Z_p(X)_{\R}}{Z_p(X)^{av}}\times
\frac{Z_p(X)^{av}}{ZQ_p(X)_{\R}}$ by $f+ZQ_p(X)_{\R} \longmapsto
(f+Z_p(X)^{av}, AP(f)+ZQ_p(X)_{\R}) $ and define $\phi:
\frac{Z_p(X)_{\R}}{Z_p(X)^{av}}\times
\frac{Z_p(X)^{av}}{ZQ_p(X)_{\R}} \longrightarrow
\frac{Z_p(X)_{\R}}{ZQ_p(X)_{\R}}$ by $(f+Z_p(X)^{av},
g+ZQ_p(X)_{\R}) \longmapsto RP(f)+g+ZQ_p(X)_{\R}$. By the two Lemmas
above, $\psi$ and $\phi$ are continuous and it is easy to check they
are inverse to each other.
\end{proof}

\section{The Generalized Harnack-Thom Theorem}
While it is easy to produce an exact sequence $H\hookrightarrow G
\longrightarrow G/H$ of topological groups, it is cumbersome to
verify that it is a locally trivial principal $H$-bundle, and worse,
it may not be in general. But the long exact homotopy sequence
induced by a fibration is extremely useful in homotopy group
calculation. We use Milnor's construction of universal bundles to
construct some weak models of the classifying spaces of some cycle
groups. They are used to produce long exact sequences of homotopy
groups. To make everything work out, we need to work in the category
of compactly generated topological spaces $\mathcal{CG}$ (see
\cite{steenrod}).

We recall that a space $X$ is compactly generated if and only if $X$
is Hausdorff and each subset $A$ of $X$ with the property that
$A\cap C$ is closed for every compact subset $C$ of $X$ is itself
closed. Since the topology of our cycle groups is defined by a
filtration of compact Hausdorff spaces, all groups we are dealing
with are in $\mathcal{CG}$. To make sure the quotient $G/H$ is in
this category, we need $H$ to be a normal closed subgroup of $G$.

Let us recall Milnor's construction of universal $G$-bundles. We
adopt the notation from page 36 of \cite{FHT}. For a topological
group $G$, let $C_G=(G\times I)/(G\times \{0\})$ be the cone on $G$,
and the $n$-th join, $G^{*n}$, is the subspace of $C_G\times \cdots
\times C_G$ of points $((g_0, t_0),..., (g_n, t_n))$ such that $\sum
t_i=1$. Thus $G^{*n}\subset G^{*(n+1)}$. For a topological group
$G\in \mathcal{CG}$, as in \cite{FHT}, we give
$E(G)=\underset{n}{\cup}G^{*n}$ the weak topology determined by
$G^{*n}$ instead of Milnor's strong topology for arbitrary
topological groups. Then we have a continuous action of $G$ in
$E(G)$ given by
$$((g_0, t_0),..., (g_n, t_n))\cdot g=((g_0g, t_0),..., (g_ng,
t_n)).$$ Set $B(G)=E(G)/G$ and let $p_G:E(G) \rightarrow B(G)$ be
the quotient map. Then
\begin{itemize}
  \item $p_G:E(G) \rightarrow B(G)$ is a principal $G$-bundle.
  \item $\pi_k(E(G))=0$ for $k\geq 0$.
\end{itemize}
The space $B(G)$ is called the classifying space of $G$ and we have
$\pi_{k+1}(B(G))\cong \pi_k(G)$. We say that a space $T$ is a weak
model of $BG$ if $T$ is weak homotopy equivalent to $BG$, i.e, they
have the same homotopy groups.

The following result is the main tool that we use to produce long
exact sequences of homotopy groups. A similar argument for
topological groups which are CW-complexes can be found in Theorem
2.4.12 of \cite{Ben}.
\begin{proposition}
Let $H, G\in \mathcal{CG}$ be two topological abelian groups and $H$
be a closed subgroup of $G$. Then we have a long exact sequence of
homotopy groups:
$$\cdots \rightarrow \pi_{k+1}(G/H)\rightarrow \pi_k(H) \rightarrow \pi_k(G)
\rightarrow \pi_k(G/H) \rightarrow \cdots$$\label{fibration}
\end{proposition}

From this result, when we have a short exact sequence $0\rightarrow
H\rightarrow G\rightarrow K\rightarrow 0$ of topological abelian
groups such that $K$ is isomorphic to $G/H$, it induces a long exact
sequence of homotopy groups. By abuse of terminology, we will call
this long exact sequence the homotopy sequence induced by the short
exact sequence.

\begin{proof}
From the principal $H$-bundle
$$\xymatrix{H\ar[r] & E(G) \ar[d]\\
& E(G)/H}$$ We see that $E(G)/H$ is a weak model of $B(H)$.

The group $G/H$, which is also in $\mathcal{CG}$ since $H$ is a
closed subgroup, acts on $E(G)/H$ and we have $(E(G/H)\times
E(G))/G=(E(G/H)\times (E(G)/H))/(G/H)$ which is a weak model of
$B(G)$. We have a fibration
$$\xymatrix{E(G)/H \ar[r] & (E(G/H)\times E(G)/H)/(G/H) \ar[d]\\
& E(G/H)/(G/H)}$$ given by the first projection which induces a long
exact sequence of homotopy groups
$$\cdots \rightarrow
\pi_{k+1}(G/H)\rightarrow \pi_k(H) \rightarrow \pi_k(G) \rightarrow
\pi_k(G/H) \rightarrow \cdots.$$ Note that since $H, G$ may not be
CW-complexes, the homomorphisms in the long exact sequence may not
be induced by maps between $H$ and $G$.
\end{proof}

We recall that for a surjective morphism $\phi:G \longrightarrow G'$
between topological groups $G, G'$, if $\phi$ is an open map or
closed map, it induces an isomorphism of topological groups $G/Ker
\phi \cong G'$. And for two normal subgroups $H_1\subset H_2 \subset
G$, we have $(G/H_1)/(H_2/H_1)\cong (G/H_2)$.

Let $A=\frac{Z_p(X)_{\R}}{2Z_p(X)_{\R}}$,
$B=\frac{Z_p(X)}{2Z_p(X)}$, $C=\frac{Z_p(X)^{av}}{2Z_p(X)_{\R}}$,
$D=\frac{Z_p(X)_{\R}}{ZQ_p(X)_{\R}}$,
$E=\frac{Z_p(X)_{\R}}{Z_p(X)^{av}}$,
$F=\frac{Z_p(X)^{av}}{ZQ_p(X)_{\R}}$,
$G=\frac{ZQ_p(X)_{\R}}{2Z_p(X)_{\R}}$.

\begin{proposition}
\begin{enumerate}
\item
The map $1+c_*:Z_p(X) \longrightarrow Z_p(X)^{av}$ is a closed map.

\item
The map $1+c_*:B \longrightarrow C$ is a closed map.

\item
The natural map $\phi:A \longrightarrow E$ is open.

\item
The natural map $\phi: A \longrightarrow D$ is open.

\item
The natural map $\phi: C \longrightarrow F$ is open.

\end{enumerate}
\label{open map}
\end{proposition}

\begin{proof}
\begin{enumerate}
\item
Let $K_1 \subset K_2 \subset \cdots \subset Z_p(X)$ be the canonical
filtration and $K^{av}_1 \subset K^{av}_2 \subset \cdots \subset
Z_p(X)^{av}$ be the canonical averaged filtration. Let $C\subset
Z_p(X)$ be a closed set. Then $(1+c_*)(C)\cap K^{av}_i=(1+c_*)(C\cap
K_{[\frac{i}{2}]})$ which is compact thus closed for all $i$.
Therefore, $(1+c_*)(C)$ is closed.

\item
Consider the following commutative diagram:
$$\xymatrix{Z_p(X) \ar[r]^{1+c_*} \ar[d]_{\pi_1} & Z_p(X)^{av} \ar[d]^{\pi_2}\\
\frac{Z_p(X)}{2Z_p(X)} \ar[r]^{1+c_*} &
\frac{Z_p(X)^{av}}{2Z_p(X)_{\R}} }$$ Let $C\subset
\frac{Z_p(X)}{2Z_p(X)}$ be a closed subset. Then $\pi_1^{-1}(C)$ is
closed which implies $\pi_2^{-1}(1+c_*)(C)=(1+c_*)(\pi_1^{-1}(C))$
is closed. Since $\pi_2$ is a quotient map, it follows that
$(1+c_*)(C)$ is closed.

\item
We have a commutative diagram
$$\xymatrix{ & \ar[ld]_{\pi_1} Z_p(X)_{\R} \ar[rd]^{\pi_2} &\\
\frac{Z_p(X)_{\R}}{2Z_p(X)_{\R}} \ar[rr]^-{\phi} & &
\frac{Z_p(X)_{\R}}{Z_p(X)^{av}}}$$

Since $\pi_1, \pi_2$ are open maps, $\phi$ is also an open map.
\end{enumerate}

The proof of (4) and (5) are similar to the proof of (3).
\end{proof}

In the following Proposition, we use the notation $T_n$ to denote
the $n$-th homotopy group of $T$ where $T$ is any of the groups $A,
..., G$. We note that all these groups are $\Z_2$-spaces so their
homotopy groups are vector spaces over $\Z_2$.

\begin{proposition}
We have the following short exact sequences:
\begin{enumerate}
\item $0\rightarrow A\longrightarrow B \longrightarrow C \rightarrow 0$

\item $0\rightarrow C\longrightarrow A \longrightarrow E \rightarrow 0$

\item $0\rightarrow G\longrightarrow A \longrightarrow D \rightarrow 0$

\item $0\rightarrow G\longrightarrow C \longrightarrow F \rightarrow 0$
\end{enumerate}

They induce long exact sequences:

\begin{enumerate}
\item $\cdots \overset{c_{n+1}}{\longrightarrow} A_n
\overset{a_n}{\longrightarrow} B_n \overset{b_n}{\longrightarrow}
C_n \overset{c_n}{\longrightarrow} A_{n-1} \longrightarrow \cdots $

\item $\cdots \overset{e_{n+1}}{\longrightarrow} C_n
\overset{c'_n}{\longrightarrow} A_n \overset{a'_n}{\longrightarrow}
E_n \overset{e_n}{\longrightarrow} C_{n-1} \longrightarrow \cdots $

\item $\cdots \overset{d_{n+1}}{\longrightarrow} G_n
\overset{g_n}{\longrightarrow} A_n \overset{a''_n}{\longrightarrow}
D_n \overset{d_n}{\longrightarrow} G_{n-1} \longrightarrow \cdots $

\item $\cdots \overset{f_{n+1}}{\longrightarrow} G_n
\overset{g'_n}{\longrightarrow} C_n \overset{c''_n}{\longrightarrow}
F_n \overset{f_n}{\longrightarrow} G_{n-1} \longrightarrow \cdots $
\end{enumerate}
\label{four long exact sequences}
\end{proposition}

\begin{proof}
(1) By Proposition \ref{smith exact sequence}, $A$ is isomorphic to
$ker(1+c_*)$. The map $1+c_*$ is surjective and closed by
Proposition \ref{open map}, thus $C$ is isomorphic as a topological
group to $B/ker(1+c_*)$. Hence we have the first exact sequence and
by Proposition \ref{fibration}, we have the first homotopy sequence.
A similar argument works for (2), (3) and (4).
\end{proof}

Since every topological abelian group is a product of Eilenberg-Mac
lane spaces, we are able to compute the homotopy types of
topological abelian groups from knowledge of their homotopy groups
alone.

\begin{example}
The homotopy types of the seven groups mentioned above for 1-cycles
on $\P^2$ are $A=K(\Z_2, 0)\times K(\Z_2, 1)\times K(\Z_2, 2),
B=K(\Z_2, 0)\times K(\Z_2,2), C=K(\Z_2, 2), D=E=K(\Z_2, 0)\times
K(\Z_2, 1), F=0, G=C$.
\end{example}

\begin{definition}
Suppose that $X$ is a real quasiprojective variety. We define the
$L_p$-total Betti number of $X$ with $\Z_2$-coefficients to be
$B(p)(X)=\sum_{k=2p}^{\infty}dim_{\Z_2} L_pH_k(X; \Z_2)$ where
$$L_pH_k(X; \Z_2)=\pi_{k-2p}(\frac{Z_p(X)}{2Z_p(X)}).$$  We define
the real $L_p$-total Betti number to be
$\beta(p)(X)=\sum_{k=p}^{\infty}dim_{\Z_2} RL_pH_k(X)$. We call
$\chi_p(X)=\sum_{k=2p}^{\infty}(-1)^k dim_{\Z_2} L_pH_k(X; \Z_2)$
the $L_p$-Euler characteristic of $X$ with $\Z_2$-coefficients and
$R\chi_p(ReX)=\sum_{k=p}^{\infty} (-1)^{k-p} dim_{\Z_2} RL_pH_k(X)$
the real $L_p$-Euler characteristic.
\end{definition}

Let $B(p)(X)_{\R}=\sum_{k=0}^{\infty}
dim_{\Z_2}\pi_k(\frac{Z_p(X)_{\R}}{2Z_p(X)_{\R}})$.

\begin{theorem}
Suppose that $X$ is a real projective variety. If $B(p)(X)$ and
$B(p)(X)_{\R}$ are finite, then
\begin{enumerate}
\item $\chi_p(X)\equiv R\chi_p(ReX) \mbox{ mod } 2 $

\item $B(p)(X)\equiv \beta(p)(X) \mbox{ mod } 2$

If in addition $G$ is weakly contractible, then

\item $B(p)(X)\geq \beta(p)(X)$.
\end{enumerate}
\label{generalized harnack-thom}
\end{theorem}

\begin{proof}
To simplify the notation, we use the same notation as in Proposition
\ref{four long exact sequences} but with different meaning. We use
$M_n$ to denote the rank of the $n$-th homotopy group of $M$,
$Kerg_n$ and $Img_n$ the rank of the kernel and the rank of the
image over $\Z_2$ of a homomorphism $g_n$ respectively.

From the finiteness assumption on $B(p)(X)$ and $B(p)(X)_{\R}$, we
know that $\sum^{\infty}_{n=0} C_n$ and $\sum^{\infty}_{n=0} E_n$
are finite from the long exact sequence 1 and 2 respectively in
Proposition \ref{four long exact sequences}.

\begin{enumerate}
\item From the first two long exact sequences in Proposition \ref{four long
exact sequences}, we have $\chi_p(B)=\chi_p(A)+\chi_p(C)$ and
$\chi_p(A)=\chi_p(C)+\chi_p(E)$, thus $\chi_p(B)\equiv \chi_p(E)
\mbox{ mod } 2$.

\item From the long exact sequence 1, we have
$A_n=Imc_{n+1}+Ima_n=C_{n+1}-kerc_{n+1}+B_n-Imb_n$ and from the long
exact sequence 2, we have
$A_n=Ima_n'+Imc'_n=E_n-kerc'_{n-1}+c_n-kerc'_n$. Simplifying the
equation and taking sums, we get $\sum B_n=\sum
E+2\sum(kerc_n-kerc'_n)$. Thus $B(p)(X)\equiv \beta(p)(X) \mbox{ mod
} 2$.

\item If $G$ is weakly contractible, then $\pi_k(G)=0$ for all $k$. From
the long exact sequence 3 and 4, we have $A_n=D_n$ and $C_n=F_n$.
Since by Theorem \ref{split}, $D_n=E_n+F_n$, we have $A_n=C_n+E_n$
for all $n$. From the long exact sequence 1, we have
$A_n=C_{n+1}-kerc_{n+1}+B_n-Imb_n$. Thus $C_n+E_n\leq C_{n+1}+B_n$.
Taking the sum over all $n$, we have $\sum E_n\leq \sum B_n$
\end{enumerate}
\end{proof}

For zero-cycles, to simplify the notation, we simply write
$B(X)=B(0)(X), B(ReX)=\beta(0)(X), \chi(X)=\chi_0(X)$ and $\chi(Re
X)=R\chi_0(Re X)$ which are the standard total Betti numbers and
Euler characteristic of $X$ and $ReX$ in $\Z_2$-coefficients.

\begin{corollary}(Harnack-Thom Theorem) Let $X$ be a real projective
variety. Then
\begin{enumerate}
\item $B(X)\geq \beta(Re X)$

\item $B(X)\equiv \beta(ReX) \mbox{ mod } 2$

\item $\chi(X)\equiv \chi(ReX) \mbox{ mod } 2 $
\end{enumerate}
\end{corollary}

\begin{proof}
For $p=0$, by Proposition \ref{o cycle},
$ZQ_0(X)_{\R}=2Z_0(X)_{\R}$, thus
$G=\frac{ZQ_0(X)_{\R}}{2Z_0(X)_{\R}}$ is trivial. By the Dold-Thom
theorem, $\pi_k(\frac{Z_0(X)}{2Z_0(X)})=H_k(X; \Z_2)$,
$\pi_k(\frac{Z_0(X)_{\R}}{2Z_0(X)_{\R}})=H_k(X/\Z_2; \Z_2)$ where
$X/\Z_2$ is the orbit space of $X$ under the action of conjugation.
Thus $B(0)(X)$ and $B(0)(X)_{\R}$ are finite. The result now follows
from the Theorem above.
\end{proof}

Suppose that $Y\subset X$ are real projective varieties and $U=X-Y$
is a real quasiprojective variety. Let $A, B, ..., G$ be the cycle
groups of $X$ defined as above and let
$A'=\frac{Z_p(Y)_{\R}}{2Z_p(Y)_{\R}}$, $B'=\frac{Z_p(Y)}{2Z_p(Y)}$,
$C'=\frac{Z_p(Y)^{av}}{2Z_p(Y)_{\R}}$,
$D'=\frac{Z_p(Y)_{\R}}{ZQ_p(Y)_{\R}}$,
$E'=\frac{Z_p(Y)_{\R}}{Z_p(Y)^{av}}$,
$F'=\frac{Z_p(Y)^{av}}{ZQ_p(Y)_{\R}}$,
$G'=\frac{ZQ_p(Y)_{\R}}{2Z_p(Y)_{\R}}$.

We can check as in Proposition 2.8 of \cite{T1} that $T'$ is
embedded as a closed subgroup of $T$ thus we identify $T'$ with its
image in $T$ for any group $T$ above. To simplify the notations, we
will use $(T')$ to mean the image of $T'$ in $T$.

\begin{lemma}
We have two short exact sequences of topological abelian groups:
\begin{enumerate}
\item $0\rightarrow \frac{A}{A'} \overset{i}{\longrightarrow} \frac{B}{B'}
\overset{1+c_*}{\longrightarrow} \frac{C}{C'} \rightarrow 0$,

\item $0\rightarrow \frac{C}{C'} \rightarrow \frac{A}{A'} \rightarrow
\frac{E}{E'} \rightarrow 0$
\end{enumerate}
and an equation
$$\frac{D}{D'}=\frac{E}{E'}\times\frac{F}{F'}.$$
\end{lemma}

\begin{proof}
We prove that the first sequence is exact and similar argument works
for the second one. Injectivity: let $a\in Z_p(X)_{\R}$ and
$a+2Z_p(X)+(\frac{Z_p(Y)}{2Z_p(Y)})\in (\frac{Z_p(Y)}{2Z_p(Y)})$,
then $a=2b+c$ where $b\in Z_p(X)$, $c\in Z_p(Y)$ and we may assume
$b, c$ have no common components and the conjugation of each
component of $c$ is not a component of $b$. Since $a$ is real, $b,
c$ are real hence $a$ is 0 in $\frac{A}{A'}$. It is trivial that the
image of $i$ is contained in the kernel of $1+c_*$ and the map
$1+c_*$ is surjective. Suppose that for $a\in Z_p(X)$,
$a+\overline{a}+2Z_p(X)_{\R}+(\frac{Z_p(Y)^{av}}{2Z_p(Y)_{\R}})\in
(\frac{Z_p(Y)^{av}}{2Z_p(Y)_{\R}})$. Since
$a+\overline{a}=2RP(a)+2AP(a)+IP(a)+IP(\overline{a})$, we have
$IP(a)+IP(\overline{a})\in Z_p(Y)^{av}$. From the definition of
$IP$, we see that $IP(a)\in Z_p(Y)$. Therefore $i$ sends
$RP(a)+AP(a)+2Z_p(X)_{\R}+(\frac{Z_p(Y)}{2Z_p(Y)})$ to
$a+2Z_p(X)_{\R}+ (\frac{Z_p(Y)^{av}}{2Z_p(Y)_{\R}})$ which implies
that the sequence is exact.
\end{proof}

By a similar calculation as in Theorem \ref{generalized
harnack-thom}, we get

\begin{theorem}
Suppose that $U$ is a real quasiprojective variety. If $B(p)(U)$ and
$B(p)(U)_{\R}$ are finite, then
\begin{enumerate}
\item $\chi_p(U)\equiv R\chi_p(ReU) \mbox{ mod } 2 $

\item $B(p)(U)\equiv \beta(p)(U) \mbox{ mod } 2$

If in addition $\frac{G}{G'}$ is weakly contractible, then

\item $B(p)(U)\geq \beta(p)(U)$.
\end{enumerate}
\label{for quasiprojective}
\end{theorem}

\begin{corollary}(Harnack-Thom Theorem for real quasiprojective varieties)
Suppose that $U$ is a real quasiprojective variety. Let $B^{BM}(U),
\beta^{BM}(ReU), \chi^{BM}(U), \chi^{BM}(ReU)$ denote the total
Betti number and the Euler characteristic of $U$ and $Re(U)$
respectively in the Borel-Moore homology with $\Z_2$-coefficients.
Then
\begin{enumerate}
\item $\beta^{BM}(ReU)\leq B^{BM}(U)$

\item $B^{BM}(U)\equiv \beta^{BM}(ReU) \mbox{ mod } 2$

\item $\chi^{BM}(U)\equiv \chi^{BM}(ReU) \mbox{ mod } 2 $
\end{enumerate}
\end{corollary}

\begin{proof}
By the Dold-Thom theorem of $\Z_2$-coefficients, we have
$$\pi_k\left(\frac{Z_0(U)}{2Z_0(U)}\right)=H^{BM}_k(U; \Z_2) \mbox{ and }
\pi_k R_0(U)=H^{BM}_k(ReU; \Z_2)$$ For zero-cycles, the equation
$\frac{D}{D'}=\frac{E}{E'}\times \frac{F}{F'}$ becomes
$\frac{A}{A'}=\frac{E}{E'}\times \frac{C}{C'}$. Now Theorem \ref{for
quasiprojective} suffices to finish the proof.
\end{proof}

\begin{example}
From the Splitting Principle in reduced real morphic cohomology (see
Theorem 5.4 in \cite{T1}) and the Duality Theorem between reduced
real morphic cohomology and reduced real Lawson homology(see Theorem
5.14 in \cite{T1}), we get
$$R_1(\P^1\times \P^1)=R_0(\P^1)\times R_1(\P^1).$$ Thus
$\pi_0(R_1(\P^1\times \P^1))=\Z_2\oplus \Z_2$, $\pi_1(R_1(\P^1\times
\P^1))=\Z_2$. $R\chi_1(\RP^1\times \RP^1)=1$, $\beta(1)(\P^1\times
\P^1)=3$. From the Splitting Principle in morphic cohomology (see
Theorem 2.10 in \cite{FL1}) and the Duality Theorem between morphic
cohomology and Lawson homology, we have
$$Z_1(\P^1\times \P^1)=Z_0(\P^1)\times Z_1(\P^1)$$ which gives
$\pi_0(\frac{Z_1(\P^1\times \P^1)}{2Z_1(\P^1\times
\P^1)})=\Z_2\oplus \Z_2$, $\pi_2(\frac{Z_1(\P^1\times
\P^1)}{2Z_1(\P^1\times \P^1)})=\Z_2$. Thus $\chi_1(\P^1\times
\P^1)=B(1)(\P^1\times \P^1)=3$.
\end{example}

\begin{example}
Consider $\C^n=\P^n-\P^{n-1}$. From the Lawson Suspension Theorem,
we are able to get $\frac{Z_p(\C^n)}{2Z_p(\C^n)}=K(\Z_2, 2(n-p)),
R_p(\C^n)=K(\Z_2, n-p)$. Therefore we have
$$L_pH_k(\C^n; \Z_2)=\left\{%
\begin{array}{ll}
    \Z_2, & \hbox{if } k=2n \\
    0, & \hbox{otherwise.} \\
\end{array}%
\right.
$$
$$RL_pH_k(\C^n)=\left\{%
\begin{array}{ll}
    \Z_2, & \hbox{if } k=n\\
    0, & \hbox{otherwise.} \\
\end{array}%
\right.$$ and then $\chi_p(\C^n)=1, R\chi_p(\R^n)=(-1)^{n-p},
B(p)(\A^n)=\beta(p)(\A^n)=1$.
\end{example}

The reduced real Lawson homology groups of a variety naturally
depend on its real structure. Two real projective varieties may be
isomorphic as complex projective varieties, but they may not be
isomorphic as real projective varieties. Thus reduced real Lawson
homology groups may be used to distinguish two real projective
varieties.

\begin{example}
Let $X$ be the smooth quadric defined by the equation
$x^2+y^2+z^2=0$ in $\P^2$. The variety $X$ is complex algebraically
isomorphic to $\P^1$ but not real algebraically isomorphic to $\P^1$
since $X$ has no real point. Therefore all the reduced real Lawson
homology groups of zero-cycles on $X$ are trivial but
$B(0)(X)=\chi_0(X)=2$.
\end{example}

\section{A construction of Weil}
Throughout this section, $X$ is a nonsingular projective variety of
dimension $m$. Let $G$ be the fundamental group of $X$ and $G'$ be
the commutator-group of $G$. Then $H=G/G'=H_1(X, \Z)$ is the first
homology group of $X$ with integral coefficients.

\begin{definition}
Let $p: (\hat{X}, y) \longrightarrow (X, x)$ be a covering map where
$\hat{X}, X$ are complex manifolds. The covering is said to be
abelian if $p_*\pi_1(\hat{X}, y)=G'$.
\end{definition}

For an abelian covering, the group of deck transformations is
isomorphic to $H$. Every element $\sigma \in H$ determines an
automorphism of $\hat{X}$, transforming each point $\hat{s}$ of
$\hat{X}$ into a point $\sigma (\hat{s})$ lying over the same point
$s$ in $X$.

Let $\epsilon$ be a group homomorphism of $H$ into the
multiplicative group $\C^*$. A multiplicative function on $\hat{X}$
with the multiplicator-set $\epsilon$ is a non-identically zero
meromorphic function $\phi$ on $\hat{X}$ such that $\phi(\sigma
\hat{s})=\phi(\hat{s})\epsilon(\sigma)$ for all $\hat{s}\in \hat{X},
\sigma\in H$.

Let $T$ be the torsion subgroup of $H$. A multiplicator-set which is
1 on $T$ is called special. Let $\Theta_X$ be the group of all
special multiplicator-sets $\epsilon$ where $|\epsilon(\sigma)|=1$
for all $\sigma\in H$. A divisor $Z$ on $X$ is defined by a
meromorphic multiplicative function $\phi$ on $\hat{X}$, as
explained in Page 873 of \cite{weil}, by taking the zero locus of
$\phi$, and the multiplicator-set of $\phi$ is special if and only
if $Z$ is algebraically equivalent to zero. It is proved in
\cite{weil} that $\Theta_X$ is an abelian variety for a nonsingular
projective variety $X$ and the real dimension of $\Theta_X$ is equal
to the rank of $H$.

\begin{definition}
For $X$ a nonsingular projective variety, the abelian variety
$\Theta_X$ is called the Picard variety of $X$.
\end{definition}

Let $Z_{m-1}(X)^{alg}$ be the group of divisors on $X$ which are
algebraically equivalent to zero.

\begin{definition}(Weil's construction)
Define a group homomorphism $w: Z_{m-1}(X)^{alg} \longrightarrow
\Theta_X$ by $$w(Z)=\epsilon$$ where $\epsilon\in \Theta_X$ is the
special multiplicator-set of $\phi$ and $Z$ is the divisor defined
by $\phi$.
\end{definition}

For the reader's convenience, we recall a definition from
\cite{weil}.

\begin{definition}
An analytic family of divisors on a nonsingular projective variety
$X$ parametrized by a nonsingular projective variety $S$ is an
algebraic cycle $V$ on $S\times X$ such that $V_s:=Pr_*(V\bullet
(s\times X))$ is a divisor on $X$ (where $Pr:S\times X
\longrightarrow X$ is the projection and $\bullet$ is the
intersection product). A mapping $f: Z_{m-1}(X)^{alg}
\longrightarrow \Theta_X$ is said to be analytic if for any analytic
family of divisors algebraically equivalent to zero on $X$,
parametrized by $S$, the map $f\circ \lambda: S \longrightarrow
\Theta_X$ is an analytic map where $\lambda: S\longrightarrow
Z_{m-1}(X)^{alg}$ is the parametrization.
\end{definition}

Let $Z_{m-1}(X)^{lin}$ be the group of divisors on $X$ which are
linearly equivalent to zero. The following is the ``Main Theorem" in
Weil's paper \cite{weil}.

\begin{theorem}
The surjective group homomorphism $w: Z_{m-1}(X)^{alg}
\longrightarrow \Theta_X$ in Weil's construction is analytic and the
kernel of $w$ is $Z_{m-1}(X)^{lin}$. There is a bijective
parametrization $\Lambda: \Theta_X \longrightarrow
\frac{Z_{m-1}(X)^{alg}}{Z_{m-1}(X)^{lin}}$. \label{weil theorem}
\end{theorem}

The main result we need is that the topology on $\Theta_X$ is
actually same as the topology on
$\frac{Z_{m-1}(X)^{alg}}{Z_{m-1}(X)^{lin}}$.

\begin{corollary}
The map $w:Z_{m-1}(X)^{alg} \longrightarrow \Theta_X$ in Weil's
construction is continuous and therefore it induces a topological
group isomorphism
$$\widetilde{w}:\frac{Z_{m-1}(X)^{alg}}{Z_{m-1}(X)^{lin}}
\longrightarrow \Theta_X.$$
\end{corollary}

\begin{proof}
We may form a topology on $Z_{m-1}(X)^{alg}$ by declaring that a set
$U\subset Z_{m-1}(X)^{alg}$ is open if and only for all
parametrizations $\lambda: S\longrightarrow Z_{m-1}(X)^{alg}$,
$\lambda^{-1}(U)$ is open. By Theorem 2.16 in \cite{Li2}, this
topology coincides with the Chow topology. Combining this with
Theorem \ref{weil theorem}, we have that $w$ is a continuous map.
From Weil's construction, $\widetilde{w} \circ \Lambda=$the identity
map. Since $\Theta_X$ is compact, $\Lambda$ is a topological group
isomorphism, which implies that $\widetilde{w}$ is a topological
group isomorphism.
\end{proof}

\begin{definition}
We say that a complex manifold $(M, -)$ with a map $-: M
\longrightarrow M$ is real if the map $-$ is an antiholomorphic map
and $-^2=$the identity map. The map $-$ is called the conjugation of
$M$.
\end{definition}

All nonsingular projective varieties defined by real polynomials in
$\P^n$ have a natural conjugation which is induced by the standard
conjugation of $\P^n$.

\begin{definition}
If $\hat{X}$ and $X$ are real complex manifolds and the covering map
$p:\hat{X} \longrightarrow X$ satisfies $p(\bar{z})=\overline{p(z)}$
for all $z\in \hat{X}$, then the covering is said to be real.
\end{definition}

\begin{lemma}
Suppose that $p:(\hat{X}, y) \longrightarrow (X, x)$ is a covering
map where $\hat{X}$ is a complex manifold and $X$ is real complex
manifold. Then the conjugation on $X$ induces a conjugation on
$\hat{X}$ such that the covering is real. \label{real covering}
\end{lemma}

\begin{proof}
For each point $t\in X$, take a small connected open neighborhood
$U_t$ and a biholomorphic local trivialization $\phi_t: p^{-1}(U_t)
\longrightarrow U_t\times F$ where $F$ is the fibre which is
discrete. We may take $U_t$ small enough and make
$U_{\bar{t}}=\overline{U_t}$ for all $t$. We define a conjugation on
$U_t\times F$ by $\overline{(w, b)}=(\overline{w}, b)$. Suppose that
$z$ is a point in the fibre over $U_t$, define
$\overline{z}=\phi^{-1}_{\bar{t}}\overline{\phi_t(z)}$. It is not
difficult to check that $\overline{\overline{z}}=z$. The map
$\phi_{\bar{t}} \circ \phi^{-1}_{\bar{s}}: \overline{U_t \cap
U_s}\times F \longrightarrow \overline{U_t \cap U_s} \times F$ is
given by $(w, b) \longrightarrow (w, gb)$ where $g$ is a bijection
from $F$ to $F$. Thus $\phi^{-1}_{\bar{t}} \circ
\phi^{-1}_{\bar{s}}$ is real, i.e., $\phi^{-1}_{\bar{t}} \circ
\phi^{-1}_{\bar{s}}\overline{(w, b)}=\overline{\phi^{-1}_{\bar{t}}
\circ \phi^{-1}_{\bar{s}}(w, b)}$. It is then easy to verify that
$\phi^{-1}_{\bar{t}}
\overline{\phi_t(z)}=\phi^{-1}_{\bar{s}}\overline{\phi_s(z)}$ in the
overlap of $U_t$ and $U_s$. So $\overline{z}$ is well defined. Since
on $U_t\times F$ for any $t$, the map sending $z$ to $\overline{z}$
is antiholomorphic, thus the map we just defined is a conjugation on
$\hat{X}$.
\end{proof}

We remark that this conjugation depends on the choice of local
trivializations, for instance in the case of trivial covering space
$\hat{X}=X\times F$.

\begin{definition}
Suppose that $p: (\hat{X}, \hat{x}) \longrightarrow (X, x)$ is a
real abelian covering and $\sigma$ is a deck transformation induced
by a loop $[f]\in \pi_1(X, x)$. Let $\gamma$ be a path from $x$ to
$\overline{x}$ and let $g={\gamma}^{-1}*\overline{f}*\gamma$ be the
loop at $x$, defined by traveling from $x$ to $\overline{x}$ along
$\gamma$, going around $\overline{x}$ along the conjugation of $f$
and then traveling back to $x$ along $\gamma$ with opposite
direction. Let $\overline{\overline{\sigma}}$ be the deck
transformation defined by $g$. If we take another path $\gamma'$
from $x$ to $\overline{x}$ and let
$g'={\gamma'}^{-1}*\overline{f}*\gamma'$, then it is easy to show
that $g'g^{-1}$ is an element in the commutator group, thus $g'$
defines a same deck transformation as $g$ does. We can check that
$$\overline{\overline{\sigma}}(z)=\overline{\sigma(\overline{z})}$$
for all $z\in \hat{X}$. We say that $\sigma$ is real if
$\sigma=\overline{\overline{\sigma}}$.
\end{definition}

From the theory of covering spaces, we know that $\sigma$ is real if
and only $[f]=[g]$ in $\pi_1(X, x)/p_*(\pi_1(\hat{X}, \hat{x}))$.

If $X$ is a nonsingular real projective manifold and $\phi$ is a
multiplicative function on $\hat{X}$, define
$\overline{\overline{\phi}}(z):=\overline{\phi(\overline{z})}$ which
is also a multiplicative function on $\hat{X}$.

Let $\epsilon\in \Theta_X$ be a special multiplicator-set and define
$\overline{\epsilon}(\sigma):=\overline{\epsilon(\overline{\overline{\sigma}})}$.
If $\epsilon=\overline{\epsilon}$, then we say that $\epsilon$ is
real. The map $\epsilon \longrightarrow \overline{\epsilon}$ induces
a conjugation on $\Theta_X$ which makes $\Theta_X$ a real complex
manifold.

\begin{lemma}
Suppose that $p:(\hat{X}, \hat{x}) \longrightarrow (X, x)$ is a real
abelian covering and $X$ is a real projective manifold. Let $Z$ be a
divisor of $X$. If $\phi$ is a multiplicative function defining $Z$,
with multiplicator-set $\epsilon$, then $\overline{\overline{\phi}}$
is a multiplicative function defining $\overline{Z}$, with
multiplicator-set $\overline{\epsilon}$. \label{real multiplicative}
\end{lemma}

\begin{proof}
$\overline{\overline{\phi}}(\sigma(y))=\overline{\phi(\overline{\sigma(y)})}=
\overline{\phi(\overline{\overline{\sigma}}(\overline{y}))}=
\overline{\phi(\overline{y})\epsilon(\overline{\overline{\sigma}})}=
\overline{\overline{\phi}}(y) \overline{\epsilon}(\sigma)$
\end{proof}

It was shown by Weil in \cite{weil} that a divisor $Z$ on $X$ is
linearly equivalent to 0 if and only if for a multiplicative
function $\phi$ defining $Z$, the special multiplicator-set
$\epsilon$ of $\phi$ is 1. Therefore, since $\overline{Z}$ is
defined by $\overline{\overline{\phi}}$ with multiplicator-set
$\overline{\epsilon}$, this implies that $\overline{Z}$ is also
linearly equivalent to 0. So the conjugation on $Z_{m-1}(X)^{alg}$
passes to $\frac{Z_{m-1}(X)^{alg}}{Z_{m-1}(X)^{lin}}$.

By Lemma \ref{real multiplicative}, it is clear that we have the
following result.

\begin{theorem}
Suppose that $p:(\hat{X}, \hat{x}) \longrightarrow (X, x)$ is a real
abelian covering and $X$ is a real projective manifold. Then the map
$\widetilde{w}: \frac{Z_{m-1}(X)^{alg}}{Z_{m-1}(X)^{lin}}
\longrightarrow \Theta_X$ in Weil's construction is real, i.e.,
$\widetilde{w}(\overline{[Z]})=\overline{w(Z)}$ for all $Z\in
Z_{m-1}(X)^{alg}$.
\end{theorem}

\begin{definition}
Suppose that $p:(\hat{X}, \hat{x}) \longrightarrow (X, x)$ is a real
abelian covering and $X$ is a real projective manifold. We say that
$X$ is real symmetric if all the deck transformations of $X$ are
real.
\end{definition}

\begin{proposition}
Suppose that $p:(\hat{X}, \hat{x}) \longrightarrow (X, x)$ is a real
abelian covering and $X$ is real symmetric. If a divisor $Z$ is
algebraically equivalent to 0, then the averaged divisor
$Z+\overline{Z}$ is linearly equivalent to 0. \label{averaged
linearly zero}
\end{proposition}

\begin{proof}
Let $\phi$ be a meromorphic multiplicative function defining $Z$,
with special multiplicator-set $\epsilon$. Since all deck
transformations are real, we have
$\overline{\epsilon}(\sigma)=\overline{\epsilon(\sigma)}$ for all
$\sigma \in H$. Then $\epsilon(\sigma)
\overline{\epsilon}(\sigma)=1$ which implies that $Z+\overline{Z}$
is linearly equivalent to 0.
\end{proof}

\begin{corollary}
A projective curve $X$ is not real symmetric if the genus $g$ of $X$
is greater than 0.
\end{corollary}

\begin{proof}
Let $p\in X$ and $D=p+\overline{p}$. For a divisor $E$ on $X$, let
$L(E)$ be the dimension of $H^0(X, [E])$ where $[E]$ is the line
bundle associated to $E$, and let $|E|$ be the linear system
associated to $E$. If $g=1$, by Riemann-Roch theorem, we have
$L(D)=2$. If $g>1$, $L(D)\geq 1$ and $L(K-D)\geq 1$ where $K$ is a
canonical divisor on $X$, then by Clifford's theorem, $L(D)\leq 2$.
Assume that for every $q\in X$, $q+\overline{q}$ is linearly
equivalently to $D$. Then $dim|D|=L(D)-1=1$. Consider the set
$\mathscr{C}_{0, 2}(X)=SP^2(X)$ of effective divisors of degree 2
where $SP^2(X)$ is the 2-fold symmetric product of $X$. We have
$\mathscr{C}_{0, 2}(X)_{\R}=SP^2(X)_{\R}=\{q+\overline{q}|q\in X\}$
and by the assumption we have $\mathscr{C}_{0, 2}(X)_{\R}\subset
|D|=\P^1$. Since the map $X\longrightarrow SP^2(X)_{\R}$ defined by
$a\longrightarrow a+\overline{a}$ is a homeomorphism, it gives an
embedding of $X$ into $\P^1$ which is impossible. Therefore, there
exists $q\in X$ such that $q+\overline{q} \notin |D|$. Since $p-q$
is algebraically equivalent to zero but
$(p-q)+\overline{(p-q)}=(p+\overline{p})-(q+\overline{q})$ is not
linearly equivalent to zero, this contradicts to the conclusion of
Corollary \ref{averaged linearly zero}. Hence, $X$ is not real
symmetric.
\end{proof}

\begin{lemma}
If $D$ is a real divisor which is linearly equivalent to 0, then
there is a real rational function $F$ such that $D=(F)$, the divisor
defined by $F$. \label{real rational}
\end{lemma}

\begin{proof}
Let $D=D_1-D_2$ where $D_1$ and $D_2$ are effective real divisors.
Since $D$ is linearly equivalent to zero, there exists a rational
functions $F=\frac{f}{g}$ such that $D=(F)$. Suppose that
$(f)=D_1+D_3$, $(g)=D_2+D_3$. Since $D_3+\overline{D}_3$ is a real
divisor, we can take a real homogeneous polynomial $h$ such that
$(h)=D_3+\overline{D}_3+D_4$. We show that we can find a real
homogeneous polynomial which defines the divisor $(fh)$. By
splitting the coefficients of $fh$ into the form $x+iy$, we may
write $fh=G+iH$ where $G, H$ are real homogeneous polynomials. If
$G$ is zero, then $(fh)=(iH)=(H)$. If $G$ is not zero, then the
degree of $G$ equals to the degree of $fh$. If $f(z)h(z)=0$, then
$f(\overline{z})h(\overline{z})=0$. We have $G(z)+iH(z)=0$ and
$G(\overline{z})+iH(\overline{z})=G(z)-iH(z)=0$. This implies that
$G(z)=0$. Thus $(G)=(fh)=D_1+2D_3+\overline{D}_3+D_4$ is a real
divisor. Similarly, we can find a real homogenous polynomial $G'$
such that $(G')=(gh)=D_2+2D_3+\overline{D}_3+D_4$. Therefore,
$D=(fh)-(gh)=(\frac{G}{G'})$.
\end{proof}

\begin{proposition}
If $D$ is a real divisor which is linearly equivalent to 0, then $D$
is in the 0-component of $Z_{m-1}(X)_{\R}$. \label{real linearly 0}
\end{proposition}

\begin{proof}
Let $(x, y)\in \C^2\backslash \{0\}$, $D=(F)$ where $F$ is a real
rational function. Let $V_{x, y}$ be the divisor defined by $x+yF$.
We have $V_{0, 1}=D$ and $V_{1, 0}=0$. Let $\gamma:[0, 1]
\longrightarrow \C^2\backslash \{0\}$ be the path given by
$\gamma(t)=(t, 1-t)$. Then each $V_{\gamma(t)}$ is real and this
gives a path in $Z_{m-1}(X)_{\R}$ joining $D$ and $0$.
\end{proof}

It follows from this result that we do not have to distinguish
between real and complex linear equivalence in $Z_{m-1}(X)_{\R}$.

Denote $Z_{m-1}(X)^0_{\R}$ to be the $0$-component of
$Z_{m-1}(X)_{\R}$ and $Z_{m-1}(X)^{lin}_{\R}=\{a \in
Z_{m-1}(X)_{\R}|a \mbox{ is linearly equivalent to } 0\}$.

\begin{corollary}
Suppose that $X$ is a real nonsingular projective variety of
dimension $m$. We have the following inclusions:
$$Z_{m-1}(X)^{lin}_{\R} \subset
Z_{m-1}(X)^0_{\R} \subset Z_{m-1}(X)^{alg}$$
\end{corollary}

\begin{definition}
Let $R_{m-1}(X)^0$ be the connected component of $R_{m-1}(X)$
containing 0 and $R_{m-1}(X)^{lin}=\{a+Z_{m-1}(X)^{av}|a \mbox{ is
linearly equivalent to some } b, \mbox{ where } b\in
Z_{m-1}(X)^{av}\}$.
\end{definition}

It is easy to check the following result.
\begin{lemma}
Suppose that $X$ is a nonsingular real projective variety of
dimension $m$. Then

$$R_{m-1}(X)^0=\frac{Z_{m-1}(X)^0_{\R}}{Z_{m-1}(X)^0_{\R}\cap
Z_{m-1}(X)^{av}},$$

$$R_{m-1}(X)^{lin}=\frac{Z_{m-1}(X)^{lin}_{\R}}{Z_{m-1}(X)^{lin}_{\R}\cap
Z_{m-1}(X)^{av}}.$$
\end{lemma}

The inclusion map $R_{m-1}(X)^{lin} \hookrightarrow R_{m-1}(X)^0$ is
a closed embedding. We will abusively denote the image of
$R_{m-1}(X)^{lin}$ in $R_{m-1}(X)^0$ by $R_{m-1}(X)^{lin}$.

Let $Pic^0(X)$ be the group of holomorphic line bundles on $X$ whose
first Chern class are zero. There is an isomorphism
$$u: \frac{Z_{m-1}(X)^{alg}}{Z_{m-1}^{lin}(X)} \longrightarrow Pic^0(X)$$
where $u$ maps a divisor $Z$ to the line bundle associated to $Z$.
We give a topology on $Pic^0(X)$ by making $u$ a homeomorphism. For
$L\in Pic^0(X)$, $L=[c]$ for some $c\in Z_{m-1}(X)^{alg}$. We define
$\overline{L}=[\overline{c}]$. Then the map $u$ is real. We have the
following commutative diagram and each map is a real topological
group isomorphism:
$$\xymatrix{\Theta_X \ar[rr]^{u \circ w} \ar[rd]^{w} &  & Pic^0(X)\\
& \frac{Z_{m-1}(X)^{alg}}{Z_{m-1}(X)^{lin}} \ar[ru]^{u} & }$$

\begin{definition}
We say that a holomorphic line bundle $L$ on a nonsingular
projective variety $X$ is real if $L$ is the line bundle associated
to some real divisor, and $L$ is averaged if $L$ is the line bundle
associated to some averaged divisor. Denote $Pic^0(X)_{\R}$ to be
the 0-component of real line bundles and $Pic^0(X)^{av}$ to be the
0-component of averaged line bundles. We define the reduced real
Picard group of $X$ to be
$$RPic^0(X)=\frac{Pic^0(X)_{\R}}{Pic^0(X)^{av}}$$ which is a topological
abelian group.
\end{definition}

The real isomorphism $u$ gives us the following result.
\begin{theorem}
For a nonsingular real projective variety $X$ of dimension $m$, we
have a topological group isomorphism:
$$\widetilde{u}: \frac{R_{m-1}(X)^0}{R_{m-1}(X)^{lin}}
\longrightarrow RPic^0(X)$$ and we get $\pi_kRPic^0(X)=0$ if $k\neq
1$, $\pi_1R_{m-1}(X)^0=\pi_1R_{m-1}(X)^{lin}\oplus \pi_1RPic^0(X)$,
and $\pi_iR_{m-1}(X)^0=\pi_iR_{m-1}(X)^{lin}$ for $i\neq 1$.
\label{rpic(X)}
\end{theorem}

\begin{proof}
Since $Pic^0(X)_{\R}, Pic^0(X)^{av}$ are closed subgroups of
$Pic^0(X)$, they are real tori. The results follow from the homotopy
sequences induced by the two exact sequences:
$$\xymatrix{ 0\ar[r] & Pic^0(X)^{av} \ar[r]& Pic^0(X)_{\R} \ar[r] & RPic^0(X) \ar[r] & 0\\
0\ar[r] & R_{m-1}(X)^{lin} \ar[r] & R_{m-1}(X)^0  \ar[r] &
RPic^0(X)\ar[r]& 0}$$
\end{proof}

We may as well as define the reduced real Picard group to be
$$RPic^0(X)=\frac{\Theta^0_{X, \R}}{\Theta^0_{X, av}}$$ where
$\Theta^0_{X, \R}$ is the 0-component of the group of
multiplicator-sets in $\Theta_X$ which are invariant under the
conjugation of $\Theta_X$, and $\Theta_{X, av}$ is the 0-component
of the group consisting of multiplicator-sets of the form $\epsilon
\overline{\epsilon}$ for $\epsilon \in \Theta_X$. From the
commutative diagram above it is easy to see that these two
topological groups are isomorphic.

\begin{corollary}
If a nonsingular real projective variety $X$ of dimension $m$ is
real symmetric, then $RPic^0(X)$ is the trivial group and therefore
$R_{m-1}^0(X)=R_{m-1}^{lin}(X)$.
\end{corollary}

\begin{proof}
Let $\epsilon:H_1(X, \Z) \longrightarrow S^1$ be a real special
multiplicator-set in $\Theta_X$. If $X$ is real symmetric, then
$\overline{\epsilon(\sigma)}=\epsilon(\sigma)$ for all $\sigma\in
H_1(X, \Z)$ which implies that $\epsilon(\sigma)=1 \mbox{ or } -1$.
Hence $\Theta_{X, \R}$ is discrete and therefore $\Theta^0_{X, \R}$
is trivial which implies that $RPic^0(X)$ is trivial and hence
$R^0_{m-1}(X)=R^{lin}_{m-1}(X)$.
\end{proof}

\section{Examples}
In this section we are going to compute the reduced real Lawson
homology groups of divisors.

Let us recall a computation done by Friedlander in \cite{F1},
Theorem 4.6.
\begin{theorem}
Suppose that $X$ is a nonsingular projective variety of dimension
$m$. Then
$$\pi_kZ_{m-1}(X)=\left\{%
\begin{array}{ll}
    NS(X), & \hbox{if } k=0; \\
    \pi_1Pic^0(X), & \hbox{if } k=1; \\
    \Z, & \hbox{if } k=2; \\
    0, & \hbox{othwerwise}.\\
\end{array}%
\right.    $$
\end{theorem}

We make a similar calculation for the real case by the method
developed in this paper.

\begin{proposition}
Suppose that $X$ is a nonsingular real projective variety of
dimension $m$. Then
$$\pi_kZ_{m-1}(X)_{\R}=\left\{%
\begin{array}{ll}
    NS(X)_{\R}, & \hbox{if } k=0; \\
    \Z_2\oplus \pi_1Pic^0(X)_{\R}, & \hbox{if } k=1; \\
    0, & \hbox{othwerwise}.\\
\end{array}%
\right.    $$ where $NS(X)_{\R}$ is the real Neron-Severi group
which is defined to be $\pi_0Z_{m-1}(X)_{\R}$. \label{real neron}
\end{proposition}

\begin{proof}
Suppose that $X\subset \P^n$. Let
$$K[X]_{\R}=\frac{\R[z_0,..., z_n]}{I_{\R}(X)}$$ be the real
coordinate ring of $X$ where $I_{\R}(X)\subset \R[z_0,..., z_n]$ is
the ideal of real polynomials vanishing over $X$. Let
$K[X]_{\R}=\oplus^{\infty}_{k=0}I_k$ where $I_k$ is the real vector
space generated by homogeneous polynomials of degree $k$ of $X$.
Define
$$K_d=\{(\frac{f}{g})|f, g \in I_d\}$$ and by Lemma \ref{real
rational}, we have a filtration
$$ \cdots \subset K_d\subset K_{d+1} \subset \cdots
=Z_{m-1}(X)^{lin}_{\R}.$$

Let $I=\lim_{d\to \infty}\P I_d\times \P I_d, \Delta_d=\mbox{
diagonal of } \P I_d \times \P I_d$, and $\Delta=\lim_{d\to \infty}
\Delta_d$ where $\P I_d$ is the real projectivisation of $I_d$, thus
$I\cong K(\Z_2\oplus \Z_2, 1)$, $\Delta\cong K(\Z_2, 1)$. For $(f_1,
g_1)\in \P I_{d_1}\times \P I_{d_1}$, $(f_2, g_2)\in \P
I_{d_2}\times \P I_{d_2}$, we define $(f_1, g_1)\cdot (f_2,
g_2):=(f_1f_2, g_1g_2)$ which induces a monoid structure on $I$ and
$\Delta$.

Let $\widetilde{I}, \widetilde{\Delta}$ be the naive group
completions of $I$ and $\Delta$ respectively. Since all $\Delta_d,
\P I_d$ are compact CW-complexes, the monoids $\Delta$ and $I$ are
free, strongly properly $c$-graded (see \cite{Li3} for the
definitions), by Theorem 4.4' of \cite{Li3}, $\widetilde{I},
\widetilde{\Delta}$ are homotopy equivalent to their homotopy
theoretic group completions respectively. Hence
$\pi_k\widetilde{I}=\pi_kI$, and $\pi_k\widetilde{\Delta}=\pi_k
\Delta$ for $k>0$ and $\pi_0\widetilde{I}=\pi_0\widetilde{\Delta}=
\Z$. Since $(I, \Delta)$ is a properly $c$-filtered free pair of
monoids, by Theorem 5.2 of \cite{Li3}, we have a fibration
$$\xymatrix{\widetilde{\Delta} \ar[r] & \widetilde{I} \ar[d]\\
& \widetilde{I}/\widetilde{\Delta}}$$ which implies that
$\widetilde{I}/\widetilde{\Delta}\cong K(\Z_2, 1)$.

There is a surjective monoid homomorphism $\phi:I \longrightarrow
Z_{m-1}(X)^{lin}_{\R}$ defined by $\phi(f, g)=(\frac{f}{g})$. We
extend it to a group homomorphism
$\widetilde{\phi}:\widetilde{I}\rightarrow Z_{m-1}(X)^{lin}_{\R}$ by
defining $\widetilde{\phi}((f, g)-(f', g'))=\phi(f, g)-\phi(f',
g')$. The kernel of $\widetilde{\phi}$ is $\widetilde{\Delta}$,
hence $Z_{m-1}(X)^{lin}_{\R}$ is isomorphic to
$\widetilde{I}/\widetilde{\Delta}$ therefore $Z_{m-1}(X)^{lin}_{\R}
\cong K(\Z_2, 1)$.

Since
$\frac{Z_{m-1}(X)^0_{\R}}{Z_{m-1}(X)^{lin}_{\R}}=Pic^0(X)_{\R}$, and
the group $Pic^0(X)_{\R}$ is a closed subgroup of $Pic^0(X)$, thus a
real torus, from the homotopy sequence induced by the short exact
sequence
$$\xymatrix{0\ar[r] & Z_{m-1}(X)^{lin}_{\R} \ar[r]&
Z_{m-1}(X)^0_{\R} \ar[r] & Pic^0(X)_{\R} \ar[r] &0}.$$

we have $$\pi_kZ_{m-1}(X)^0_{\R}=\left\{%
\begin{array}{ll}
    \Z_2\oplus \pi_1Pic^0(X)_{\R}, & \hbox{if } k=1; \\
    0, & \hbox{otherwise}. \\
\end{array}%
\right.$$

This completes the proof.
\end{proof}

\begin{proposition}
For a nonsingular real projective variety $X$, $NS(X)_{\R}$ is
finitely generated.
\end{proposition}

\begin{proof}
Let $m$ be the dimension of $X$ and let
$$H=\frac{Z_{m-1}(X)_{\R}\cap
Z_{m-1}(X)^{alg}}{Z_{m-1}(X)^{lin}_{\R}}.$$ By Proposition \ref{real
linearly 0}, $H$ is embedded as a closed subgroup of
$Pic^0(X)=\frac{Z_{m-1}(X)^{alg}}{Z_{m-1}(X)^{lin}}$ hence $H$ is a
Lie subgroup of $Pic^0(X)$ therefore $\pi_0(H)$ is finitely
generated. The inclusion map $i:Pic(X)_{\R} \rightarrow Pic(X)$
induces a map $i_*:NS(X)_{\R}=\pi_0(Pic(X)_{\R})\rightarrow
\pi_0(Pic(X))$ whose kernel is $\pi_0(H)$, and because
$\pi_0(Pic(X))$ is finitely generated hence $NS(X)_{\R}$ is finitely
generated.
\end{proof}

Let
$$T=\{a-b|a+\overline{a}=b+\overline{b}, a, b\in
\mathscr{C}_{m-1}(X)\}.$$

\begin{lemma}
$T=\{a-\overline{a}|a\in \mathscr{C}_{m-1}(X)\}$.
\end{lemma}

\begin{proof}
Suppose that $a, b \in \mathscr{C}_{m-1}(X)$ and
$a+\overline{a}=b+\overline{b}$. Write $a=\sum^n_{i=1}n_iV_i$ where
each $V_i$ is an irreducible subvariety and $n_i>0$. Since $a-b\in
T$, we may assume that $a$ and $b$ have no common irreducible
subvariety components. From the relation
$a+\overline{a}=b+\overline{b}$, we see that each $V_i$ must be a
component of $\overline{b}$. Thus $b=\overline{a}$.
\end{proof}

The following observation is the main tool that we are going to use
to compute $R_{m-1}(X)$.

\begin{proposition}
We have the following exact sequences of topological groups:
\begin{enumerate}
\item $0\longrightarrow Z_{m-1}(X)_{\R} \hookrightarrow
Z_{m-1}(X)\overset{Sa}{\longrightarrow} T \longrightarrow 0$

\item $0\longrightarrow T \hookrightarrow Z_{m-1}(X)
\overset{Av}{\longrightarrow} Z_{m-1}(X)^{av} \longrightarrow 0 $
\end{enumerate}
where $Sa(c)=c-\overline{c}$ and $Av(c)=c+\overline{c}$, the groups
$T$ and $Z_{m-1}(X)^{av}$ are isomorphic as a topological group to
$\frac{Z_{m-1}(X)}{Z_{m-1}(X)_{\R}}$ and $\frac{Z_{m-1}(X)}{T}$
respectively. \label{two exact seq}
\end{proposition}

\begin{proof}
A direct verification shows that the sequences are exact. To show
that $T$ is isomorphic as a topological group to
$\frac{Z_{m-1}(X)}{Z_{m-1}(X)_{\R}}$, it suffices to prove that $Sa$
is a closed map. Let $K_1\subset K_2 \subset \cdots \subset
Z_{m-1}(X)$ be the canonical filtration. The topology of $T$ is the
subspace topology of $Z_{m-1}(X)$. For $C\subset Z_{m-1}(X)$ a
closed subset, $C\cap K_n$ is compact and $Sa(C\cap K_n)=Sa(C)\cap
K_{2n}$ which is closed for any $n$, so $Sa$ is a closed map. The
map $Av$ is a closed map which is proved in Proposition \ref{open
map}.
\end{proof}

\begin{lemma}
Suppose that $X$ is a nonsingular real projective variety of
dimension $m$.
\begin{enumerate}
\item
For $k\geq 2$, $\pi_k Z_{m-1}(X)^{av}$ is a 2-torsion group.

\item
The homotopy group
$$\pi_kT=\left\{
          \begin{array}{ll}
            0, & \hbox{ if } k > 2; \\
            \Z, & \hbox{ if } k=2; \\
            \pi_1(\frac{T}{Z_{m-1}(X)^{lin}}), & \hbox{ if } k=1\\
          \end{array}
        \right.
$$ and hence $\pi_kT$ is free for $k>0$.

\item
If $H^1(X, \C)=0$ and $NS(X)$ is free, then $\pi_0T$ is free.
\end{enumerate}
\label{2-torsion}
\end{lemma}

\begin{proof}
\begin{enumerate}
\item
For a continuous map $f:S^k\longrightarrow Z_{m-1}(X)^{av}$ where $k
\geq 2$, $2f(t)\in 2Z_{m-1}(X)_{\R} \subset Z_{m-1}(X)^{av}$. But
since $\pi_k(2Z_{m-1}(X)_{\R})\cong \pi_k(Z_{m-1}(X)_{\R})=0$ for
$k\geq 2$, we see that $2f$ is null homotopic. Thus
$\pi_kZ_{m-1}(X)^{av}$ is a 2-torsion group for $k\geq 2$.

\item
Let $T^0$ be the zero-component of $T$. Since $T$ is a closed
subgroup of $Z_{m-1}(X)$,
$\frac{T^0}{Z_{m-1}(X)^{lin}}\hookrightarrow
\frac{Z_{m-1}(X)^{alg}}{Z_{m-1}(X)^{lin}}=Pic^0(X)$ is a closed
embedding and hence $\frac{T^0}{Z_{m-1}(X)^{lin}}$ is a closed Lie
subgroup of $Pic^0(X)$ which implies that
$\pi_k\frac{T}{Z_{m-1}(X)^{lin}}$ is free for $k>0$. By a similar
calculation of the homotopy type of $Z_{m-1}(X)^{lin}_{\R}$ in
Proposition \ref{real neron}, we get $Z_{m-1}(X)^{lin}\cong K(\Z,
2)$. From the homotopy sequence induced by the short exact sequence
$$\xymatrix{0 \ar[r] & Z_{m-1}(X)^{lin} \ar[r] & T \ar[r]&\frac{T}{Z_{m-1}(X)^{lin}}\ar[r] & 0}$$
we get the result.

\item
Since $H^1(X, \C)=0$, $Pic^0(X)$ is trivial, hence
$Z_{m-1}(X)^{alg}=Z_{m-1}(X)^{lin}$. Then
$\pi_0Z_{m-1}(X)=\frac{Z_{m-1}(X)}{Z_{m-1}(X)^{lin}}$ is free from
the hypothesis. From the homotopy sequence induced by the short
exact sequence above, we see that $\pi_0T$ is also free.
\end{enumerate}
\end{proof}

\begin{theorem}
For a nonsingular real projective $X$ of dimension $m$,
$\pi_kR_{m-1}(X)=0$ for $k>2$. \label{vanish}
\end{theorem}

\begin{proof}
Consider the two exact sequences in Proposition \ref{two exact seq}
\begin{enumerate}
\item $0\longrightarrow Z_{m-1}(X)_{\R} \hookrightarrow
Z_{m-1}(X)\overset{Sa}{\longrightarrow} T \longrightarrow 0$

\item $0\longrightarrow T \hookrightarrow Z_{m-1}(X)
\overset{Av}{\longrightarrow} Z_{m-1}(X)^{av} \longrightarrow 0 $
\end{enumerate}

From the homotopy sequence induced by the first exact sequence, we
see that $\pi_kT=0$ if $k>2$. And from the homotopy sequence induced
by the second exact sequence, we get $0\rightarrow \pi_3
Z_{m-1}(X)^{av} \overset{\phi_1}{\rightarrow} \pi_2T
\overset{i_{2*}}{\rightarrow} \pi_2Z_{m-1}(X)\rightarrow
\pi_2Z_{m-1}(X)^{av} \overset{\phi_2}{\rightarrow} \pi_1T
\rightarrow \pi_1Z_{m-1}(X) \rightarrow \cdots$. By Lemma
\ref{2-torsion}, the groups $\pi_1T, \pi_2T$ are free and the groups
$\pi_2Z_{m-1}(X)^{av}$, $\pi_3Z_{m-1}(X)^{av}$ are 2-torsion groups,
thus $\phi_1, \phi_2$ are 0-maps and $\pi_3Z_{m-1}(X)^{av}=0$.

Consider the composition of the maps: $ T \overset{i_2}{\rightarrow}
Z_{m-1}(X) $ $\overset{\rho_1}{\rightarrow} T$ which is $2 id:
T\longrightarrow T$, and consider the induced maps on homotopy
groups: $\pi_2T \overset{i_{2*}}{\longrightarrow}
\pi_2Z_{m-1}(X)\overset{\rho_{1*}}{\longrightarrow} \pi_2T$. All of
these groups are isomorphic to $\Z$ and $\rho_{1*}\circ
i_{2*}=2id_*$ which is of degree 2. If $\rho_{1*}:
\pi_2Z_{m-1}(X)\longrightarrow \pi_2T$ is of degree 1, from the
homotopy sequence induced by the first exact sequence, we have
$0\rightarrow \pi_2Z_{m-1}(X) \overset{\rho_{1*}}{\rightarrow}
\pi_2T \rightarrow \pi_1Z_{m-1}(X)_{\R} \rightarrow \pi_1Z_{m-1}(X)
\rightarrow \cdots $ which implies that $\pi_1Z_{m-1}(X)_{\R}$ is
mapped injectively into $\pi_1Z_{m-1}(X)$ which is impossible. Thus
$\rho_{1*}$ is of degree 2 and $i_{2*}$ is of degree 1. Since
$\phi_2$ is a 0-map, we have $\pi_2Z_{m-1}(X)^{av}=0$.

The result now follows from the homotopy sequence induced by
$$\xymatrix{0\ar[r]& Z_{m-1}(X)^{av} \ar[r] & Z_{m-1}(X)_{\R} \ar[r]& R_{m-1}(X)\ar[r] & 0}$$
\end{proof}

\begin{definition}
Suppose that $X$ is a nonsingular real projective variety. The
Picard number $\rho(X)$ is the rank of the free part of
$\pi_0Z_{m-1}(X)$, and we call $RNS(X)=\pi_0R_{m-1}(X)$ the reduced
real Neron-Severi group of $X$.
\end{definition}

\begin{corollary}
Suppose that $X$ is a nonsingular real projective variety of
dimension $m$, we have the following relation:
$$\rho(X)+1\equiv
dim_{\Z_2}RNS(X)+dim_{\Z_2}\pi_1R_{m-1}(X)+dim_{\Z_2}\pi_2R_{m-1}(X)
mod 2$$
\end{corollary}

\begin{proof}
By the fundamental theorem on finitely generated abelian groups,
$\pi_0Z_{m-1}(X)=\Z^k\oplus \Z^t_2\oplus G$ where $k=\rho(X)$, $t$
is a nonnegative integer, and $G$ is in the torsion part of
$\pi_0Z_{m-1}(X)$ where $\Z_2$ is not a direct summand of $G$.

From the homotopy sequence induced by
$$\xymatrix{0\ar[r] & 2Z_{m-1}(X) \ar[r] & Z_{m-1}(X) \ar[r]
& \frac{Z_{m-1}(X)}{2Z_{m-1}(X)}\ar[r] & 0 }$$ and the homotopy
groups of $Z_{m-1}(X)$, we have an exact sequence: $0\rightarrow
\pi_22Z_{m-1}(X) \rightarrow \pi_2Z_{m-1}(X) \rightarrow
\pi_2(\frac{Z_{m-1}(X)}{2Z_{m-1}(X)}) \rightarrow \pi_12Z_{m-1}(X)
\rightarrow \pi_1Z_{m-1}(X) \rightarrow
\pi_1(\frac{Z_{m-1}(X)}{2Z_{m-1}(X)}) \rightarrow \pi_02Z_{m-1}(X)
\rightarrow \pi_0Z_{m-1}(X) \rightarrow
\pi_0(\frac{Z_{m-1}(X)}{2Z_{m-1}(X)})\rightarrow 0$. The map $i_*:
\pi_*2Z_{m-1}(X) \rightarrow \pi_*Z_{m-1}(X)$ is easily seen to be a
map of degree 2, hence
$\pi_2(\frac{Z_{m-1}(X)}{{2Z_{m-1}(X)}})=\Z_2$,
$dim_{\Z_2}\pi_1(\frac{Z_{m-1}(X)}{2Z_{m-1}(X)})=rkPic^0(X)+t$, and
$dim_{\Z_2}\pi_0(\frac{Z_{m-1}(X)}{2Z_{m-1}(X)})=k+t$. Therefore
$\sum_{n=0}^2dim_{\Z_2}\pi_n(\frac{Z_{m-1}(X)}{2Z_{m-1}(X)})=1+rkPic^0(X)+2t+k$
and we note that the rank of the Picard variety $Pic^0(X)$ is even.
The result now follows from Theorem \ref{vanish} and Theorem
\ref{generalized harnack-thom}.

\end{proof}

\begin{corollary}
Suppose that $X$ is a nonsingular projective variety of dimension
$m$. We assume that $H^1(X, \C)=0$ and $NS(X)$ is free. Then we have
$$\rho(X)+1\equiv dim_{\Z_2}RNS(X)+dim_{\Z_2}\pi_1R_{m-1}(X)mod 2.$$
\label{last}
\end{corollary}

\begin{proof}
By Lemma \ref{2-torsion}, $\pi_0T$ is free. Since $Pic^0(X)$ is
trivial, $\pi_1Z_{m-1}(X)=0$. Consider the homotopy sequence induced
by the second exact sequence of Proposition \ref{two exact seq}, we
have $$\cdots \rightarrow \pi_1Z_{m-1}(X) \rightarrow
\pi_1Z_{m-1}(X)^{av} \rightarrow \pi_0T \rightarrow 0$$ which
implies that $\pi_1Z_{m-1}(X)^{av}$ is free. From the homotopy
sequence induced by the exact sequence $0\rightarrow Z_{m-1}(X)^{av}
\rightarrow Z_{m-1}(X)_{\R} \rightarrow R_{m-1}(X)\rightarrow 0$ and
Theorem \ref{vanish}, we have an exact sequence $0\rightarrow
\pi_2R_{m-1}(X) \rightarrow \pi_1Z_{m-1}(X)^{av}$, but since $\pi_2
R_{m-1}(X)$ is 2-torsion, this implies $\pi_2R_{m-1}(X)=0$.
\end{proof}

\begin{corollary}
If $X$ is a real complete intersection of dimension $>2$, then
$$dim_{\Z_2}RNS(X)\equiv dim_{\Z_2}\pi_1R_{m-1}(X) mod 2.$$
\end{corollary}

\begin{proof}
By the weak Lefschetz theorem and the exponential sequence on $X$,
we have $NS(X)\cong \Z$ and $H^1(X, \C)=0$, so $\rho(X)=1$. The
result then follows from Corollary \ref{last}.
\end{proof}

\begin{acknowledgement}
The author thanks Blaine Lawson for his encouragement and patience
in listening to the proof and Christian Haesemeyer for some helpful
remarks. He is also indebted to the referee for careful
proofreading, corrections and suggestions, the National Center for
Theoretical Sciences of Taiwan at Hsinchu for providing a wonderful
working environment, and Yusuf Mustopa for help with some linguistic
matters.
\end{acknowledgement}

\bibliographystyle{amsplain}

\end{document}